   \documentclass[12pt,onecolumn,twoside]{IEEEtran}
   
\IEEEoverridecommandlockouts                              

\usepackage{graphicx}

\usepackage{url}
\usepackage{amsmath,amssymb,mathtools}
 \usepackage{amsthm,thmtools}
\usepackage{amsfonts} 

\usepackage{enumerate}
\usepackage{color}             \usepackage{lipsum}     

\usepackage{verbatim}

\usepackage{amssymb}
\usepackage{amsbsy}
\usepackage{amsmath}
\usepackage{setspace}
\usepackage{url}

\usepackage{subcaption}
\usepackage{float}
 \usepackage{dsfont}
 
\usepackage{tikz}
\usetikzlibrary{positioning}

 \newcommand{\diag}{\operatorname{diag}}

\newcommand{\trace}{\operatorname{trace}}

 \declaretheorem[qed=$\blacksquare$ ]{Example}

\declaretheorem[name={Theorem}  ] {Theorem}

\declaretheorem[name={Remark}  ] {Remark}
\declaretheorem[name={Corollary}  ] {Corollary}

\declaretheorem[name={Proposition}  ] {Proposition}

\newcommand {\R}{\mathbb R}

\newcommand{\be}{\begin{equation}}
\newcommand{\ee}{\end{equation}}

\DeclareMathOperator{\vol}{vol}

\usepackage{lineno}

 \title{\LARGE \bf
A sufficient condition for $k$-contraction\\ of the series connection  of two systems }

\author{Ron Ofir,   Michael Margaliot, Yoash Levron, and
 Jean-Jacques Slotine
\thanks{This research was  partially supported by a research grant from the Israel Science Foundation~(ISF).}
\thanks{RO and YL are  with 
		the Andrew and Erna Viterbi Faculty of Electrical Eng., Technion---Israel Institute of Technology, Haifa 3200003, Israel.}%
\thanks{MM (Corresponding Author) is  with the School of Electrical  Eng.,
		and the Sagol School of Neuroscience, 
		Tel-Aviv University, Tel-Aviv~69978, Israel.
		E-mail: \texttt{michaelm@tauex.tau.ac.il}}%
		\thanks{JJS is with the Department of Mechanical Eng. and
        the Department of Brain and Cognitive Sciences, Massachusetts Institute of Technology, Cambridge, Massachusetts, USA.}
}
\date{\today}

\begin{document}

\maketitle

\begin{abstract}
    The flow  of contracting   systems  contracts 1-dimensional parallelotopes, i.e., line segments, at an exponential rate.  One reason for the  usefulness of contracting  systems  is that many interconnections of contracting sub-systems yield an overall  contracting system. 
    
    A  generalization of contracting systems  is~$k$-contracting systems, where~$k\in\{1,\dots,n\}$. 
    The flow  of such systems contracts the volume of $k$-dimensional parallelotopes  at an exponential rate, and in particular they reduce to contracting systems when~$k=1$. It was shown by
    Muldowney and Li that time-invariant  $2$-contracting systems have a well-ordered asymptotic  behaviour: all bounded trajectories converge to the  set of equilibria.     
    
    Here, we   derive a sufficient condition guaranteeing that the system obtained from the series interconnection of two sub-systems is~$k$-contracting.  This is based on a new formula for the $k$th multiplicative and additive compounds of a block-diagonal matrix, which
    may be of independent interest. As an application, we find conditions guaranteeing that $2$-contracting systems with an exponentially
    decaying input   retain the well-ordered behaviour of  time-invariant  2-contracting systems. 
\end{abstract}


\section{Introduction}
Contraction theory plays an important role in systems and control theory, with applications in robotics, systems biology, neuroscience, the design of observers and controllers, and more. This is due to several reasons. First, contraction implies a highly-ordered asymptotic behavior:  any two trajectories approach one another at an exponential rate~\cite{LOHMILLER1998683}. Therefore, if an equilibrium point exists it is unique and globally exponentially stable. In a time-varying and $T$-periodic contracting system all trajectories converge to a unique and globally asymptotically stable limit cycle with period~$T$, and the convergence occurs at an exponential rate~\cite{LOHMILLER1998683,entrain2011,RFM_entrain}. Thus, the system \emph{entrains} to the periodic excitation modeled by the $T$-periodic vector field. In fact, nonlinear contracting
systems have a well-defined frequency response, as shown in~\cite{freq_convergent} in the context of convergent  systems~\cite{pavlov_book}. 
Second, there  exist  sufficient conditions for contraction based on matrix measures~\cite{LOHMILLER1998683,sontag_contraction_tutorial}. These conditions are easy to verify, and just like Lyapunov's second theorem,    do not require solving the differential equations. However, 
there is a growing awareness  of the importance of finding the ``correct metric''  for establishing 
contraction~\cite{LOHMILLER1998683, Mier_coop_finsler,forni2014,sontag2010contractive}.
Third, various interconnections of contracting sub-systems, including parallel, series, and some feedback connections,  yield an overall contracting system~\cite{LOHMILLER1998683,network_contractive}. 

Contraction theory is an active area of research, with recent contributions including both theoretical    results and applications. 
Recent generalizations of contraction theory  typically guarantee   weaker  conclusions on the asymptotic behaviour, but allow studying a larger class of systems. Some extensions include various notions of ``weak contraction'' (see, e.g.~\cite{jafarpour2020weak,cast_book}), contraction of piecewise-smooth  dynamical systems~\cite{contraction_non_diff}, and the introduction of~$\alpha$-contracting systems~\cite{pines2021}, with~$\alpha\geq 1$ \emph{real}, which is motivated in part by the seminal
works of Douady and   Oesterl\'{e}~\cite{Douady1980}, and 
Leonov  and  his   colleagues (see the recent monograph by Kuznetsov and   Reitmann~\cite{book_volker2021}) who developed powerful tools to  bound  the Hausdorff dimension of complex attractors in chaotic dynamical systems. 

One generalization, which is the focus of this work, is called $k$-contraction, and was suggested in~\cite{kordercont} (see also the note~\cite{weak_manchester}). Let~$n$ denote the dimension of the dynamical  system, and fix~$k\in\{1,\dots,n\}$. The flow of a~$k$-contracting system contracts the volume of~$k$-dimensional  parallelotopes at an exponential rate.  In particular, for~$k=1$ these are just standard contracting  systems. Ref.~\cite{kordercont} also provides   simple to check sufficient conditions for~$k$-contraction based on the $k$th additive compound  of the Jacobian  of the vector field. 

The notion of~$k$-contraction is motivated in part by the seminal work of
  Muldowney and his colleagues~\cite{muldo1990,li1995}, on systems that, using the new terminology, are~$2$-contracting. One of the  main results in this context is that every bounded solution of a time-invariant  $2$-contracting system converges to an equilibrium point. This is strictly weaker than the case of $1$-contracting systems,
  as the equilibrium point is not necessarily unique. Also,~$2$-contracting systems with a  time-varying $T$-periodic vector fields do not necessarily entrain. The theory of~$2$-contracting systems  proved useful in the analysis of epidemiological models (see, e.g.~\cite{SEIR_LI_MULD1995}). Indeed, 
  such models typically have  at least two  equilibrium points, corresponding
  to the disease-free and the endemic steady-states, and thus are not~$1$-contracting with respect to (w.r.t.) any norm. 
  
Since various  interconnections of contracting sub-systems yield an overall contracting system, it is natural to ask if the same holds for~$k$-contracting systems as well~\cite{weak_manchester}.
Ref.~\cite{ron_CDC2021} considered the series interconnection of~$k$-contracting systems with~$k\in\{1,2\}$, and showed that such interconnections are in general not $2$-contracting.
Furthermore, the asymptotic analysis of such systems is more delicate than in the case of $1$-contracting systems because the well-ordered behaviour of $2$-contracting systems
only holds in the time-invariant case, while connecting two sub-systems implies that at least one sub-system has an input from the other sub-system and is thus time-varying.

Here, we derive a new sufficient condition guaranteeing that the series interconnection of two sub-systems is~$k$-contracting for some~$k\in\{1,\dots,n\}$. Roughly speaking, this condition requires  that both sub-systems are~$k$-contracting and, in addition, they must satisfy  ``additive $j$-contraction conditions'', with~$j<k$. We show using an example that in general this condition cannot be improved.

Our main result  is based on a new closed-form expression for the~$k$th multiplicative compound and the~$k$th additive compound of a block-diagonal matrix. We believe that this result may find other applications as well. 
Matrix compounds have recently found many applications in systems and control theory, see the tutorial paper~\cite{comp_barshalom}. 

The remainder of this   paper is organized as follows. The next section reviews known definitions and results that are used later on. Section~\ref{sec:main}
includes the main theoretical results.
Section~\ref{sec:appli} describes an application of these results, and the final section concludes.  

We use standard notation. Small [capital]  letters denote column vectors [matrices]. $I_n$ is the~$n\times n$ identity matrix, and we drop the subscript if the dimension is clear from the context. 
For a matrix~$A$,~$A^T$ is the transpose of~$A$.
For a square matrix~$A$,~$\det(A)$
 is the determinant   of~$A$, and $\trace(A)$ is the trace of~$A$.

\section{Preliminaries} 
The sufficient  condition for~$k$-contraction
    in~\cite{kordercont}
is based on the~$k$th  additive compound of the Jacobian of the vector field. 
To make this paper  more accessible, we briefly review matrix compounds, their geometrical interpretation,  and their role in the analysis of ODEs. 
For more details, see also~\cite{muldo1990}. For more  recent applications of these compounds in systems and control theory, see~\cite{margaliot2019revisiting,rami_osci,pines2021,cheng_diag_stab} and the tutorial paper~\cite{comp_barshalom}.

\subsection{Multiplicative compound}
Let~$C\in\R^{n\times m}$. Fix~$k\in\{1,\dots,\min\{n,m\}\}$. 
Recall that a $k$ minor of~$C$ is the determinant of some~$k\times k$ submatrix of~$C$. 
The   \emph{$k$th multiplicative compound} of~$C$, denoted~$C^{(k)}$, is the~$\binom{n}{k}\times\binom{m}{k}$ matrix that contains all the~$k\times k$ minors of~$C$ in lexicographic order~\cite{muldo1990}. For example, for~$n=m=3$ and~$k=2$,
$C^{(2)}$ is the~$3\times 3$ matrix:
\[
C^{(2)}= \begin{bmatrix}
    \det(\begin{smallmatrix}  c_{11}& c_{12} \\ c_{21}& c_{22} \end{smallmatrix} ) &
         \det(\begin{smallmatrix}  c_{11}& c_{13} \\ c_{21}& c_{23} \end{smallmatrix})  &
             \det(\begin{smallmatrix}  c_{12}& c_{13} \\ c_{22}& c_{23} \end{smallmatrix})  \\
    \det( \begin{smallmatrix}  c_{11}& c_{12} \\ c_{31}& c_{32} \end{smallmatrix}  )&
        \det( \begin{smallmatrix}  c_{11}& c_{13} \\ c_{ 31}& c_{33} \end{smallmatrix})  &
           \det(  \begin{smallmatrix}  c_{12}& c_{13} \\ c_{32}& c_{33} \end{smallmatrix})  \\
  \det( \begin{smallmatrix}  c_{21}& c_{22} \\ c_{31}& c_{32} \end{smallmatrix} ) &
         \det(\begin{smallmatrix}  c_{21}& c_{23} \\ c_{31}& c_{33} \end{smallmatrix} ) &
            \det( \begin{smallmatrix}  c_{22}& c_{23} \\ c_{32}& c_{33} \end{smallmatrix}  )
\end{bmatrix}.
\]
This definition implies that $(C^{(k)})^T = (C^T)^{(k)}$, $C^{(1)}=C$, and if~$n=m$ then~$C^{(n)}=\det(C)$. Also, $(I_n)^{(k)} = I_r$, where~$r := \binom{n}{k}$.

The Cauchy–Binet  formula~\cite[Chapter~0]{matrx_ana},  asserts  that for any~$B \in\R^{n\times m},C\in\R^{m\times p}$ and any~$k\in \{1,\dots, \min \{n,m,p\}\}$,
we have
\be\label{eq:mcp}
(BC)^{(k)}=B^{(k)}  C^{(k )}.
\ee
This justifies the term multiplicative compound. 
In particular,~\eqref{eq:mcp} 
implies that if~$n=m=p$ then~$\det(BC)=\det(B)\det(C)$, and that if~$A\in \R^{n\times n}$ is non-singular then
\[
    (I_n)^{(k)} = (A A^{-1})^{(k)} = A^{(k)} (A^{-1})^{(k)},
\]
so~$(A^{(k)})^{-1}=(A^{-1})^{(k)}$.

One  reason for the importance of the multiplicative compound is that it has a useful geometric  interpretation. 
Fix~$k\in\{1,\dots,n\}$ and vectors~$x^1,\dots,x^k \in \R^n$. The parallelotope
generated by these vectors (and the zero vertex) is
\[
P(x^1,\dots,x^k):=\{ \sum_{i=1}^k r_i x^i: r_i \in [0,1] \}.
\]
Note that this implies that~$0\in\R^n$ is a vertex of~$P$. 
This  can always be assured  by a simple translation.
Collect the vectors in the matrix~$X:=\begin{bmatrix}
x^1&\dots & x^k
\end{bmatrix}\in\R^{n\times k}$. The  \emph{Gram matrix} associated with~$x^1,\dots,x^k$
is the $k\times k$ symmetric matrix:
\begin{align}\label{eq:defgram}
    G ( x^1,\dots,x^k): &=X^T X \nonumber\\
    &=
    \begin{bmatrix}
    (x^1) ^T x^1 & (x^1) ^Tx^2 & \hdots & (x^1) ^T x^k \\
     &   \vdots \\
     (x^k) ^T x^1 & (x^k) ^Tx^2 & \hdots & (x^k) ^T x^k 
    \end{bmatrix}.
\end{align}
It follows from~\eqref{eq:defgram}
  that for any~$r\in\R^k$ we have
$
|\sum_{i=1}^k r_i x^i  |^2 = r^T G(x^1,\dots,x^k) r,
$
so~$G(x^1,\dots,x^k)$ is non-negative definite, 
and   it is positive-definite if and only if~$x^1,\dots,x^k$ are linearly independent. 
The 
  volume of the parallelotope~$P(x^1,\dots,x^k)$ 
 is
\be\label{eq:volpeqn}
\vol(P(x^1,\dots,x^k))=\sqrt{\det(G(x^1,\dots,x^k))},
\ee
(see, e.g.~\cite[Chapter~IX]{Gantmacher_vol1}). 
For example, suppose that~$x^i=e^i$, $i=1,\dots,k$, where~$e^i$ is the~$i$th canonical vector in~$\R^n$. Then~\eqref{eq:volpeqn} gives
\begin{align*}
(\vol(P(e^1, & \dots,e^k)))^2 \\&= \det(G(e^1,\dots,e^k))\\
&= \det \left( \begin{bmatrix}
(e^1)^T e^1 & (e^1)^T e^2 & \dots & (e^1)^T e^k \\
& \vdots & \\
(e^k)^T e^1 & (e^k)^T e^2 & \dots & (e^k)^T e^k 
\end{bmatrix} \right) \\
&= 1. 
\end{align*}

The volume of~$P(x^1,\dots,x^k)$
can be expressed  
using the $k$th multiplicative compound of~$X$. To see this, note that  combining~\eqref{eq:defgram} and the Cauchy-Binet formula yields
 \begin{align*}
     \det(G(x^1,\dots,x^k))&=
     \det(X^TX)\\
     &=(X^T X) ^{(k)}\\
     &=(X^T) ^{(k)} X ^{(k)}\\
     &=(X ^{(k)})^T X ^{(k)}.
 \end{align*}
Since~$X\in\R^{n\times k}$, 
 $X^{(k)}$ is an~$\binom{n}{k}$ column vector, so~$\det(G(x^1,\dots,x^k))= |X ^{(k)}|^2$.
 Thus, 
\be\label{eq:vol_via_compound}
\vol(P(x^1,\dots,x^k))=|X ^{(k)}|. 
 \ee
 We emphasize again that here~$X ^{(k)}$ is a column vector. The norm~$|\cdot|$ here is the~$L_2$ norm, but using the equivalence of norms in~$\R^n$ other norms can be used when studying asymptotic 
 properties like convergence to zero of the volume of a parallelotope when the vectors~$x^1,\dots,x^k$ evolve in time.  
 
In the special case where~$k=n$, $X$ is a square matrix, and~\eqref{eq:vol_via_compound}
reduces to  the well-known formula
\begin{align*}
\vol(P(x^1,\dots,x^k))      &=|\det(X))|\\&=
|\det(\begin{bmatrix}
    x^1&\dots & x^n
    \end{bmatrix})|.
\end{align*}

The multiplicative compound has a  useful  spectral property. Let~$A\in\R^{n\times n}$ with 
eigenvalues~$\lambda_i$, $i\in\{1,\dots,n\}$. The  eigenvalues of~$A^{(k)}$ are all the products
$
   \lambda_ {i_1} \lambda_{i_2} \dots \lambda_{i_k  } ,
$
with~$1\leq i_1<i_2 <\dots <i_k\leq n   $.
The next example demonstrates this. 
\begin{Example}
Consider the matrix
\[
A=\begin{bmatrix}
a_{11} & a_{12} & a_{13} & a_{14} \\
0 & a_{22} & a_{23} & a_{24} \\
0&0& a_{33} & a_{34} \\
0&0&0 & a_{44}
\end{bmatrix}.
\]
A calculation gives that~$A^{(3)}$ is the matrix depicted in Fig.~\ref{fig:a3}
\begin{figure*}
\[
A^{(3)}=\begin{bmatrix}
a_{11} a_{22} a_{33}& a_{11} a_{22} a_{34}& a_{11} (a_{23} a_{34}-a_{24} a_{33} )&
 a_{14} a_{22} a_{33} - a_{12} a_{24}
 a_{33} - a_{13} a_{22} a_{34} + a_{12} a_{23} a_{34}\\
 0&
 a_{11} a_{22} a_{44}& a_{11} a_{23} a_{44}& a_{12} a_{23} a_{44}-a_{13} a_{22} a_{44} \\
 0& 0& 
 a_{11} a_{33} a_{44}& a_{12} a_{33} a_{44}\\
 0& 0& 0& a_{22} a_{33} a_{44}
\end{bmatrix} 
\]
\caption{The matrix~$A^{(3)}$.}\label{fig:a3}
 \hrule
 \end{figure*}
and thus the eigenvalues of~$A^{(3)} $ are the products of three eigenvalues of~$A$. 
\end{Example}

Suppose now that the vectors~$x^1,\dots,x^k$ vary with time via~$\dot x^i(t) = A x^i(t)$, and let 
\begin{align*}
X(t) &:= \begin{bmatrix} x^1(t) & \dots & x^k(t)\end{bmatrix} \\&= \exp(At)\begin{bmatrix} x^1(0) & \dots & x^k(0)\end{bmatrix}.
\end{align*}
Then the norm of~$X^{(k)}(t) = (\exp(At))^{(k)} \begin{bmatrix} x^1(0) & \dots & x^k(0) \end{bmatrix}^{(k)}$ is the volume of the parallelotope generated by these vectors at time~$t$. This naturally leads to the following question: how does~$(\exp(At))^{(k)}$ evolve in time? To address this, we require the~$k$th additive compound of~$A$.

\subsection{Additive compound}
Let~$A\in\R^{n\times n}$. 
For~$k\in\{1,\dots,n\}$, 
the   \emph{$k$th additive compound} of~$A $    is the~$\binom{n}{k}\times\binom{n}{k}$ matrix defined by
\[
A^{[k]}:=\frac{d}{d t } (\exp(At) )^{(k)} |_{t=0}.
\]
This implies that~$(I+ t 
 A)^{(k)}=I+t A^{[k]} +o(t)$. In particular,
 $
 A^{[1]}=A $, 
  and  $A^{[n]}=\trace(A).
 $
 
 It can be shown using the definition of the additive compound and the properties of the multiplicative compound that if~$A,B\in\R^{n\times n}$ then
 \[
 (A+B)^{[k]}=A^{[k]}+B^{[k]}.
 \]
 This justifies the term additive compound.
 
If~$A$ has eigenvalues~$\lambda_1,\dots,\lambda_n$ then the eigenvalues of~$A^{[k]}$  are all the sums 
$
   \lambda_ {i_1}+ \lambda_{i_2} +\dots+ \lambda_{i_k  } ,
$
with~$1\leq i_1<i_2 <\dots <i_k\leq n$~\cite{muldo1990}.

A useful relation between the multiplicative and additive compounds under the matrix exponential is~\cite{muldo1990}
\begin{equation}\label{eq:mul_add_exp}
    \exp(A)^{(k)} = \exp(A^{[k]}).
\end{equation}

\subsection{Compounds and ODEs}
In the context of dynamical systems, the importance of these compounds  is due to following fact. If~$\Phi:\R _+\to \R^{n\times n}$ is the solution of the matrix differential equation
\[
\frac{d}{dt} \Phi(t)=A(t) \Phi(t),\quad \Phi(0)=I, 
\]
where~$t \to A(t)$ is continuous, then
\be\label{eq:rs}
\frac{d}{dt} (\Phi(t))^{(k)}=(A(t))^{[k]}
(\Phi(t))^{(k)}. 
\ee
In other words,~$(\Phi(t))^{(k)}$ also evolves according to a linear dynamics, with the matrix~$(A(t))^{[k]}$.
Roughly speaking,~$A^{[k]}$ determines the evolution of $k$-dimensional parallelotopes under the LTI dynamics~$\dot x=Ax$~\cite{wini2010}.

A time-varying nonlinear system is called $k$-contracting 
if its variational equation along any solution (which is an LTV) contracts $k$-dimensional parallelograms;
see~\cite{kordercont} for the exact definition. For our purposes, it is enough to review a sufficient condition for $k$-contraction. 
Recall that a vector norm~$|\cdot|:\R^n\to\R_+$ induces a matrix norm~$||A||:=\max_{|x|=1} |Ax|$, and a matrix measure~$\mu(A):=\lim_{\varepsilon\to 0^+}  (||I+\varepsilon A||-1)/ \varepsilon$. If~$\mu( (A(t))^{[k]} )\leq-\eta<0$ all~$t\geq 0$ then
applying Coppel's inequality~\cite{coppel1965stability} to~\eqref{eq:rs} yields
 $|| (\Phi(t))^{(k)} ||\leq \exp(-\eta t) 
|| (\Phi(0))^{(k)} ||$ for all~$t\geq 0$. This leads to the following. 
\begin{Proposition}\cite{kordercont}
Consider the time-varying nonlinear system~$\dot x(t)=f(t,x(t))$, with~$f$ a~$C^1$ mapping,
and suppose that its trajectories evolve on a convex set~$\Omega\subseteq \R^n$. Let~$J(t,x):=\frac{\partial}{\partial x} f(t,x)$ denote the Jacobian of~$f$ with respect to~$x$. If
\be\label{eq:gpo}
\mu \left( (J(t,z))^{[k]} \right )\leq -\eta<0,\text{ for all }
t\geq 0,z\in \Omega.
\ee
then the system is~\emph{$k$-contracting}.
\end{Proposition}
Note that for~$k=1$ this reduces to the  standard infinitesimal contraction  condition~\cite{sontag_contraction_tutorial},
as~$J^{[1]}=J$. Thus,
a~$1$-contracting system is a contracting 
system.

Note also that condition~\eqref{eq:gpo} is robust in the sense that if it holds for~$f$ then it also holds for small perturbations of~$f$ (but perhaps with a different value~$\eta)$.

For~$p\in\{1,2,\infty\}$, let~$\mu_p$   denote the matrix measure induced by the~$L_p$
vector norm~$|\cdot|_p$. An important advantage of contraction theory  is that that there exist easy to verify sufficient conditions for contraction in terms of matrix measures. 
It is useful to provide similar conditions for $k$-contraction. 
These  can be easily derived using    the following result. Let~$Q(k,n)$ denote the set of all increasing sequences of~$k$ numbers from~$\{1,\dots,n\}$ ordered lexicographically. For example,
\[
Q(3,4)=\{
\{1,2,3\},\{1,2,4\},\{1,3,4\},
\{2,3,4\}
\}.
\]

\begin{Proposition} (see, e.g.~\cite{muldo1990})
Let~$A\in\R^{n\times n}$. Fix~$k\in\{1,\dots,n \}$.
Then the 
 $L_1$, $L_2$, and~$L_\infty$
 matrix measures of~$A^{[k]}$ are:
\small
\begin{align}\label{eq:matirxm_k}
\mu_1(A^{[k]}) &= 
 \max_{\alpha\in Q(k,n)
 }\left ( \sum_{p=1}^k a_{\alpha_p,\alpha_p}   + \sum_{\substack{j \notin \alpha }}(|a_{j,\alpha_1}| + \cdots + |a_{j,\alpha_k}|)\right ) ,\nonumber \\
\mu_2(A^{[k]}) &= \sum_{i=1}^k \lambda_i\left( \frac{A + A^T}{2} \right) , \\
\mu_{\infty}(A^{[k]}) &=
\max_{\alpha \in Q(k,n)}\left(  \sum_{p=1}^k a_{\alpha_p,\alpha_p} + \sum_{\substack{j \notin\alpha  }}(|a_{\alpha_1,j}| + \cdots + |a_{\alpha_k,j}|)\right) , \nonumber
\end{align}
\normalsize where for a symmetric matrix~$S\in\R^{n\times n}$,  
$
\lambda_1(S)\geq\dots\geq\lambda_n(S) , 
$
are the eigenvalues of~$S$.
\end{Proposition}

\subsection{Contraction with a hierarchic norm}
To analyze contraction for block-diagonal matrices, we use an interesting result of Str{\"o}m~\cite{strom1975} (see also~\cite{network_contractive}). For the sake of completeness, we briefly review this result. Given~$x\in\R^s$, decompose it as 
\begin{equation}
    x = \begin{bmatrix}
        x^1 \\
        \vdots \\
        x^r 
    \end{bmatrix},
\end{equation}
where $x^i \in \R^{s_i}$, and~$ \sum_{i=1}^r s_i =s $. Let $|\cdot|_i$ denote a norm on $\R^{s_i}$, and let~$|\cdot|_0$ denote a monotonic norm on~$\R^r$ (see~\cite{mono_norms} for the definition and properties of monotonic norms).
Define a norm~$|\cdot|:\R^ s  \to \R_+$ by
\begin{equation}\label{eq:hier_norm}
    |x| := \left| \begin{array}{c}
        |x^1|_1 \\
        \vdots \\
        |x^r|_r
    \end{array} \right|_0.
\end{equation}
This may be interpreted as a ``hierarchic norm'', as we first decompose~$x$ into sub-vectors and then combine the norms of all these sub-vectors to form $|x|$.

Given $B \in \R^{s \times s}$, partition it  into $r \times r$ blocks $B^{ij} \in \R^{s_i \times s_j}$, and define their matrix norms by
\begin{equation}
    ||B^{ij}||_{ij} := \sup_{  z \in \R^{s_j} \setminus\{0\} } \frac{|B^{ij} z|_i}{|z|_j}.
\end{equation}

\begin{Theorem}{\cite[Thm.~9]{strom1975}} \label{thm:strom}
Let~$\mu$ denote the matrix measure induced by the norm $|\cdot|$ defined in~\eqref{eq:hier_norm}. 
Let~$\mu_i$ denote the matrix measure induced by $|\cdot|_i$, $i = 0,\dots,r$.
 Define $C \in \R^{r \times r}$ by
    \begin{equation}
        c_{ij}: = \begin{cases}
            \mu_i(B^{ii}), & i = j, \\
            ||B^{ij}||_{ij}, & i \neq j.
        \end{cases}
    \end{equation}
    Then
    \begin{equation}\label{eq:lub}
        \max_i \mu_i(B^{ii}) \le \mu(B) \le \mu_0(C).
    \end{equation}
\end{Theorem}
This provides in particular an upper bound on~$\mu(B)$, for the ``big'' matrix~$B\in\R^{s\times s}$, using~$\mu_0(C)$, with the smaller matrix~$C\in\R^{r\times r}$. If $\mu_0(C) \le -\eta < 0$ then~\eqref{eq:lub} implies that $\dot{x} = Bx$ is contracting.

Note that if $B$ is block-diagonal, that is,~$B^{ij}=0$ for all~$i\not = j $, then $C$ is diagonal and
the lower and upper bounds in~\eqref{eq:lub}
are equal, since for any matrix measure~$\mu$ induced by a monotonic norm and any square diagonal matrix~$E$,
\begin{equation}
    \mu(E) = \max_i \lambda_i,  
\end{equation}
where~$\lambda_i$ is the~$i$th eigenvalue of~$E$. 
Therefore, for block-diagonal matrices Thm.~\ref{thm:strom} yields 
\begin{equation}\label{eq:strom_diag}
    \mu(B) = \max_i \mu_i(B^{ii}).
\end{equation}

      The next section describes our main results.
\section{Main Results} \label{sec:main}

Consider the series interconnection of two time-varying  nonlinear sub-systems   
\begin{align}\label{eq:oversys}
\dot x^1&=f^1(t,x^1),\nonumber\\
\dot x^2&=f^2(t,x^2, x^1),
\end{align}
with~$x^1\in\R^n$ and~$x^2\in \R^m$.
We assume that the trajectories of this system evolve in a convex state-space~$\Omega^1\times\Omega^2 \subseteq \R^n\times \R^m$.
We also assume that~$f^1,f^2$ are~$C^1$.  The Jacobian of~\eqref{eq:oversys}
is
 \be\label{eq:jaco_serial}
 J(t,x)=\begin{bmatrix}J^{11} (t,x^1) & 0 \\
        J^{21}(t,x) & J^{22}(t,x) \end{bmatrix},
 \ee
where~$x:=\begin{bmatrix}x^1\\x^2\end{bmatrix}$ and~$J^{ij} := \frac{\partial}{\partial x^j} f^i$. From here on, we always assume that the initial time is~$t=0$. 

\subsection{Sufficient condition for $k$-contraction of  a series interconnection}

We can now state our first main result. For a square matrix~$M$, we define the~$0$th compounds by~$M^{(0)} := 1$ and~$M^{[0]} := 0$.
\begin{Theorem}\label{thm:serial_k_contract}
    Assume that~$J^{21}(t,x )$ is uniformly bounded for all~$x  \in \Omega^1 \times \Omega^2$ and~$t \ge 0$. Fix~$k \in \{1,\dots,n+m\}$, and let $i_1 := \max\{0,k-n\}$,~$i_2 := \min\{m,k\}$.
    If there exist matrix measures $\mu_i$ induced by (possibly different) $L_p$ norms and constants $\eta_i > 0$ such that
    \begin{equation}\label{cond:block_contraction}
        \mu_i((J^{11}(t,x^1))^{[k-i]}) + \mu_i((J^{22}(t,x))^{[i]}) \le -\eta_i < 0.
    \end{equation}
    for all $i \in \{i_1,\dots,i_2\}$, $x  \in \Omega^1\times  \Omega^2$ and $t \ge 0$,
    then~\eqref{eq:oversys} is $k$-contracting.
\end{Theorem}

\begin{Remark}
For~$k=1$, $i_1=0$ and~$i_2=1$, so~\eqref{cond:block_contraction} becomes
\begin{align*}
    \mu((J^{11}(t,x^1))^{[1]}) &\le -\eta_1 , \\
    \mu((J^{22}(t,x))^{[1]}) &\le -\eta_2,
\end{align*}
(where for simplicity we use the same matrix measure). Thus, the sufficient condition for contraction is that both sub-systems are contracting, This is a well-known result     (see for example~\cite{LOHMILLER1998683} or~\cite{sontag2010contractive}). 
  
For~$k=2$ (and assuming that~$n,m\geq 2$), $i_1=0$ and~$i_2=2$, so~\eqref{cond:block_contraction} becomes
\begin{align}\label{eq:casek=2}
    \mu((J^{11}(t,x^1))^{[2]}) &\le -\eta_1 , \nonumber\\
\mu((J^{11}(t,x^1))^{[1]}) + \mu((J^{22}(t,x))^{[1]})&\le -\eta_2 ,\\
    \mu((J^{22}(t,x))^{[2]}) &\le -\eta_3. \nonumber
\end{align}
Thus, the sufficient condition for $2$-contraction of the series interconnection~\eqref{eq:oversys} is that both sub-systems are~$2$-contracting,
and also the ``additive $1$-contraction condition''~\eqref{eq:casek=2}. Note that this condition implies that for any~$(t,x)$ at least one of the sub-systems is~$1$-contracting. 
\end{Remark}

A fundamental  issue   in systems theory is
the complexity vs. stability problem: ``Does an increase of complexity lead to an improvement of
system stability, or is it the other way around?'' (see, e.g.~\cite{may_nature_1972,siljak_lrge_scale}). Thm.~\ref{thm:serial_k_contract}
 shows that $k$-contraction of the series connection requires more than $k$-contraction of each sub-system.  In this respect, 
  an increase in the 
  complexity of the system   can 
  only destroy  the $k$-contraction property, for any~$k>1$.
  
\begin{Example}
Consider the series connection
\begin{align}\label{eq:lti:ser}
\dot x^1&=Ax^1, \nonumber \\
\dot x^2&=B x^1+C x^2,
\end{align}
with~$x^1,x^2\in \R^2$,
$A=\diag(1,-2)$, $B=0$,
and~$C=\diag(\zeta_1,\zeta_2)$, with~$\zeta_1\geq \zeta_2$ and~$\zeta_1+\zeta_2<0$. Both sub-systems are $2$-contracting, as~$\trace(A)=-1$ and $\trace(C)=\zeta_1+\zeta_2$. Let~$e^i$ denote the $i$th canonical vector in~$\R^4$, and consider the 2D time-varying parallelotope
\[
P(x(t, e^1),x(t,e^3) )=\{   r_1 x(t, e^1) + r_3x(t, e^3): r_i \in [0,1] \}.
\]
Since~$x(t,e^1)=\exp(t)e^1$ and~$x(t,e^3)= \exp(\zeta_1 t) e^3 $,
the volume of~$P(x(t,e^1),x(t,e^3))$ evolves like~$\exp((1+\zeta_1)t )$. Thus, a necessary condition for~$2$-contraction is that~$1+\zeta_1<0$
and this is exactly the 
additive 
  $1$-contraction condition.
This shows that the sufficient  condition in Theorem~\ref{thm:serial_k_contract} cannot be improved in the general case.
\end{Example}

\begin{Remark}
It is important to note however that  
 the additive 
  $1$-contraction condition in~\eqref{eq:casek=2}
  does not require  that either 
  of the two sub-systems is contracting on the entire state-space.  
To illustrate this, consider the series interconnection
\begin{equation}
\begin{aligned}
    \dot{x}_1 &= -\frac{1}{2}x_1^2 - x_1, \\
    \dot{x}_2 &= x_2 x_1,
\end{aligned}
\end{equation}
with $\Omega = \R^2_+$. The 
Jacobian of this system is
$
J(x)=\begin{bmatrix}
-x_1 -1 & 0\\
x_2 & x_1
\end{bmatrix},
$
and~$J^{[2]}=\trace(J)=-1$, so the system is~$2$-contracting.
Yet,  the dynamics of  $x_1$   is not contracting on the entire state-space, and the same holds for the dynamics of~$x_2$.
In this case, the 
 additive 
  $1$-contraction condition is
  \[
  \mu( -x_1-1)+\mu(x_1 ) \leq-\eta<0,
  \]
  and for~$\eta=1$ this indeed holds on the entire state-space. 
\end{Remark}

Note that for  the special case  $k=2$ and  the matrix measure~$\mu_2$ induced by the $L_2$ norm, Theorem~\ref{thm:serial_k_contract} reduces to~\cite[Thm.~3]{weak_manchester}.
  
The remainder of this section is devoted to the proof of Theorem~\ref{thm:serial_k_contract}. This is based on an auxiliary result that describes the $k$th compounds of a block-diagonal matrix.

\subsection{ Compounds of a block-diagonal matrix}
For $A_i \in \R^{n_i \times n_i}$, $i=1,\dots,s$,
let 
\[
\diag(A_1,\dots,A_s):=\begin{bmatrix}
    A_1 &0 &0& \dots & 0 \\
    0&A_2 & 0& \dots&0 \\
    & &\vdots\\
    0&0&\dots &0& A_s
\end{bmatrix}.
\] 
  We use $\otimes$  to denote
the Kronecker product.
The Kronecker sum~\cite[Ch. 4]{HornJohnsonTopics} of~$A\in\R^{n\times n}$ and~$B\in\R^{m\times m}$ is
\begin{equation}
    A \oplus B := A \otimes I_m + I_n \otimes B.
\end{equation}

We can now state our second main result. 
This 
  relates the multiplicative [additive] compound matrix of block-diagonal matrices to Kronecker products [sums] of the compound matrices of the individual blocks.

\begin{Theorem}\label{thm:ron}
Let~$A \in \R^{n \times n}$ and $B \in \R^{m \times m}$.   Define~$C:=\diag(A,B)$. Pick~$k \in \{1,\dots,n+m\}$, and let~$r:=\binom{n+m}{k}$.
Let~$i_1:=\max\{0,k-n\}$, and~$i_2:=\min\{m,k\}$.
There exists a permutation matrix~$P\in\R^{r\times r}$ such that
    \begin{align}\label{eq:blk_diag_mul}
          C^{(k)}   
        &= P \left( \diag_{i \in \{ i_1,\dots,i_2 \}}(A^{(k-i)} \otimes B^{(i)}) \right) P^{-1},
    \end{align}
    and
    \begin{align}\label{eq:blk_diag_add}
         C^{[k]} &= P \left( \diag_{i \in \{ i_1,\dots,i_2 \}}(A^{[k-i]} \oplus B^{[i]}) \right) P^{-1},
    \end{align}
    where
    for a  square matrix~$M$, we define $M^{(0)} := 1$ and $M^{[0]} := 0$.
\end{Theorem}
 
\begin{Remark}
It is useful to verify that the dimensions of the matrices in~\eqref{eq:blk_diag_mul} agree. Clearly, $C^{(k)}$ has dimensions~$\binom{n+m}{k}\times \binom{n+m}{k}$.
The matrix on the right-hand side in~\eqref{eq:blk_diag_mul} has dimensions~$\ell\times\ell$, where
\begin{align*}
    \ell&:= \sum_{i = \max\{0,k-n\}}^{\min\{m,k\}} \binom{n}{k-i}\binom{m}{i} \\
    &= \sum_{i = 0}^{k} \binom{n}{k-i}\binom{m}{i}\\
    &= \binom{n+m}{k},
\end{align*}
where the second equality holds since any potentially added term is zero, and the third equality is Vandermonde's identity.
\end{Remark}

\begin{Example}
Note that for~$k=1$ we have~$i_1=0$ and~$i_2=1$ so both~\eqref{eq:blk_diag_mul} and~\eqref{eq:blk_diag_add} give
$C  =P  \diag( A , B) P^{-1}$ which indeed holds for~$P=I$, whereas for~$k=n+m$ we get
\begin{align*}
  C^{(n+m)} & = P  \left (A ^{(n)}\otimes B^{(m)}\right ) P^{-1}\\
 &=P\det(A)\det(B)P^{-1}\\
 &= \det(A)\det(B) ,
\end{align*} 
as the only~$1\times 1$ permutation matrix is the scalar one, and similarly
\begin{align*}
  C^{[n+m]}  &=P \left ( A ^{[n]} \oplus B^{[m]}\right )  P^{-1}\\
                  &=  \trace(A) + \trace(B) .
\end{align*}
\end{Example}

\begin{Example}\label{exa:diagcase}
Suppose that~$n=3$, $A=\begin{bmatrix}
    \lambda_1 &0 &0\\
    0&\lambda_2&0\\
    0&0&\lambda_3
\end{bmatrix}$,
$m=2$ and~$B=\begin{bmatrix}
    \lambda_4 &0 \\
    0&\lambda_5
\end{bmatrix}$.
Let~$C:=\diag(A,B)$ and take~$k=2$. A direct  calculation yields that~$C^{(2)}$ is diagonal, with diagonal entries
\begin{align*}
    \lambda_1 \lambda_2,\lambda_1 \lambda_3,\lambda_1 \lambda_4,\lambda_1 \lambda_5,\lambda_2 \lambda_3,\lambda_2 \lambda_4,\lambda_2 \lambda_5, 
    \lambda_3 \lambda_4, \lambda_3 \lambda_5,\lambda_4 \lambda_5.
\end{align*}
On the other-hand,
\[
\diag( A^{(2)}, A^{(1)}\otimes B^{(1)}, B^{(2)} )
\]
is a diagonal matrix with diagonal 
entries
\begin{align*}
    \lambda_1 \lambda_2,\lambda_1 \lambda_3,\lambda_2 \lambda_3,
    \lambda_1 \lambda_4,
    \lambda_1 \lambda_5,
    \lambda_2 \lambda_4,
    \lambda_2 \lambda_5, 
    \lambda_3 \lambda_4, 
    \lambda_3 \lambda_5,
    \lambda_4 \lambda_5.
\end{align*}
This implies that there exits a permutation matrix~$P\in\R^{10\times 10}$ such that~\eqref{eq:blk_diag_mul} holds.
\end{Example}

\begin{IEEEproof}[Proof of Theorem~\ref{thm:ron}]
We
begin by proving~\eqref{eq:blk_diag_mul}.
Let~$Q(k,n,m)$ denote
the set of increasing sequences of~$k$ integers in $\{1, \dots,n+m\}$. (We write~$Q(k,n,m)$ rather than just~$Q(k,n+m)$ because we define below a  
 special lexicographic ordering of the elements in~$Q(k,n,m)$
 that depends on the parameter~$n$). 
 The cardinality of~$Q(k,n,m)$ is~$\binom{n+m}{k}$. For a sequence $\alpha=\{\alpha_1,\dots,\alpha_k\} \in Q(k,n,m)$,
let $s_\alpha$ denote the minimal $i \in \{1,\dots,k\}$ such that $\alpha_i > n$, or $k+1$ if no such~$i$ exists. 
In other words,~$s_\alpha<k+1$ implies that~$\alpha$ includes at least one index~$\alpha_\ell>n$ and since~$A\in\R^{n\times n}$, 
this index corresponds to  a row in~$C$ that is a row from~$B$. Note that~$i_i+1 \le s_\alpha \le i_2+1$.

We define a ``block-lexicographic ordering''    by  $\alpha \in Q(k,n,m)$ precedes $\beta  \in Q(k,n,m)$ if
\begin{enumerate}
    \item $s_\alpha > s_\beta$, or
    \item $s_\alpha = s_\beta$, and $\alpha$ precedes $\beta$ in the standard lexicographic ordering.
\end{enumerate}
For example, if~$n=3$, $m=2$, and~$k=2$ then the   ordering is
\begin{align*}
    Q(2,3,2)=\{ & \{1,2\},\{1,3\},\{2,3\},\{1,4\},\\
    &
\{1,5\}, \{2,4\},\{2,5\},\{3,4\},\{3,5\},\{4,5\}
\},
\end{align*}
whereas, for~$n=2$, $m=3$, and~$k=2$  the   ordering is
\begin{align*}
Q(2,2,3)=\{ & \{1,2\},\{1,3\},\{1,4\},
\{1,5\},\\& \{2,3\},\{2,4\},\{2,5\},\{3,4\},\{3,5\},\{4,5\}
\}, 
\end{align*}
which in this case   
is just the lexicographic ordering on~$Q(2,5)$.

Let~$D^{(k)}:=P^{-1} C^{(k)} P$, where~$P \in \R^{r \times r}$ is the permutation matrix such that the entries of~$D^{(k)}$ are those of~$C^{(k)}$, but now 
organized in the block-lexicographic ordering.
We first show that each entry of~$D^{(k)}$ is the determinant of a block-triangular matrix, or equivalently that for any~$\alpha, \beta \in Q(k,n,m)$, the submatrix~$C[\alpha|\beta]$ is a block-triangular matrix. Fix~$\alpha,\beta \in Q(k,n,m)$ such that $s_\alpha \ge s_\beta$. Decompose the~$k \times k$ matrix~$C[\alpha|\beta]$ as the block matrix
\begin{equation}\label{eq:C_blk}
    C[\alpha|\beta] = \begin{bmatrix}
        C^{11} & C^{12} \\
        C^{21} & C^{22}
    \end{bmatrix},
\end{equation}
where
\begin{align*} 
C^{11} & := C[\{\alpha_1,\dots,\alpha_{s_\alpha-1}\}|\{\beta_1,\dots,\beta_{s_\alpha-1}\}] , \\
C^{12} & := C[\{\alpha_1,\dots,\alpha_{s_\alpha-1}\}|\{\beta_{s_\alpha}, \dots, \beta_{k}\}]
 .
\end{align*}
Note that
$C^{11}   \in  
\R^{(s_\alpha-1 )\times (s_\alpha-1)}$,
and
$C^{12}  
\in  
\R^{(s_\alpha-1 ) \times (k-s_\alpha+1)}
$.

By definition, $s_\alpha$ is the minimal index~$i$ such that~$\alpha_i>n$, so the block~$C^{11}$ is a submatrix of~$A$ with $(s_\alpha - s_\beta)$ additional columns of zeros added to the right, and the block~$C^{22}$ is a submatrix of~$B$. Furthermore, every entry of the block~$C^{12}$ is zero, since~$\alpha_{s_\alpha - 1} \le n$ and~$\beta_{s_\alpha} > n$.

   We consider three cases. First, assume that~$s_\alpha = 1$, then~$C[\alpha|\beta] = C^{22}$ and~$C(\alpha|\beta) = \det(C^{22})$. Second, assume that~$s_\alpha = k+1$, then~$C[\alpha|\beta] = C^{11}$ and~$C(\alpha|\beta) = \det(C^{11})$. Finally, assume that $1 < s_\alpha < k+1$. In this case, all blocks are defined, but since all entries of~$C^{12}$ are zero, we have that~$C(\alpha|\beta) = \det(C^{11})\det(C^{22})$.

We now show that~$D^{(k)}$ has a block-diagonal structure. This is equivalent to showing that~$C (\alpha|\beta)=0$ whenever~$s_\alpha>s_\beta$ or~$s_\alpha<s_\beta$. Consider the case~$s_\alpha > s_\beta$, then the block~$C^{11}$ has at least one column of zeros, so~$C(\alpha|\beta) = 0$. The proof that~$C(\alpha|\beta) = 0$ when~$s_\alpha < s_\beta$ is similar.

Consider the diagonal blocks of~$D^{(k)}$. 
The  set
\[
\{ D  (\alpha|\beta) : s_\alpha=s_\beta = \ell\} 
\]
is   the \emph{$\ell$th diagonal block} of~$D^{(k)}$, and its dimensions are $\left( \binom{n}{\ell-1} \binom{m}{k-\ell+1} \right) \times \left( \binom{n}{\ell-1} \binom{m}{k-\ell+1} \right)$. Consider an entry of the $\ell$th diagonal block for some $\ell \in \{i_1+1,\dots,i_2+1 \}$, and fix~$\alpha,\beta$
such that~$s_\alpha = s_\beta = \ell$. Since~$s_\alpha = s_\beta$, the block~$C^{11}$ is a submatrix of~$A$ and $C(\alpha|\beta) = \det(C^{11}) \det(C^{22})$ is a product of an~$(\ell - 1)$ minor of~$A$ and a~$(k - \ell + 1)$ minor of~$B$. Therefore,~$C(\alpha|\beta)$ is an entry of~$A^{(\ell - 1)} \otimes B^{(k - \ell + 1)}$.

We now show that the~$\ell$th diagonal block of~$D^{(k)}$ is just~$A^{(\ell-1)} \otimes B^{(k-(\ell-1))}$.
Let $\alpha \in Q(k,n,m)$ such that $s_\alpha = \ell$. Then $\alpha = \{\alpha^A, \alpha^B\}$ where $\alpha^A = \{\alpha_1, \dots, \alpha_{\ell-1}\}$ is a sequence of indices ``pointing'' to entries of $A$, and~$\alpha^B = \{\alpha_\ell, \dots, \alpha_k\}$ is a sequence of indices ``pointing'' to entries of $B$. The sequence $\beta = \{\beta^A, \beta^B\} \in Q(k,n,m)$ that immediately succeeds $\alpha$ according to the block-lexicographic ordering may be found as follows: if $\alpha^B$ is not the last element in $Q(k-\ell+1,m)$ according to the lexicographic ordering, i.e. $\alpha^B \neq \{m-k+\ell-1, \dots, m\}$, choose $\beta^B$ to be the element following $\alpha^B$ in $Q(k-\ell+1,m)$ and $\beta^A = \alpha^A$. Otherwise, choose $\beta^B$ to be the first element in $Q(k-\ell+1,m)$ and $\beta^A$ to be the element succeeding $\alpha^A$ in $Q(\ell,n)$. Following this, since the cardinality of $Q(k-\ell+1,m)$ is $\binom{m}{k-\ell+1}$, we conclude that if $\alpha$ is the $i$th element of $Q(k,n,m)$ in the block-lexicographic ordering such that $s_\alpha = \ell$, then $\alpha_A$ is element number $\lceil i / \binom{m}{k-\ell+1} \rceil$ in $Q(\ell,n)$ and $\alpha^B$ is element number $(i-1) \% \binom{m}{k-\ell+1} + 1$ in $Q(k-\ell+1,m)$, where $a \% b$ denotes the remainder of the integer division $a / b$. Recalling that
\begin{equation}
    (A \otimes B)_{i,j} = (A)_{\lceil i / m \rceil, \lceil j / m \rceil} (B)_{(i-1) \% m + 1, (j-1) \% m + 1},
\end{equation}
we conclude that the $\{i,j\}$ entry of the $\ell$th diagonal block of $D^{(k)}$ is equal to the $\{i,j\}$ entry of~$A^{(\ell-1)} \otimes B^{(k-(\ell-1))}$.
This completes the proof of~\eqref{eq:blk_diag_mul}.


To prove~\eqref{eq:blk_diag_add},
let~$H_i:=\exp(At)^{(k-i)} \otimes \exp(Bt)^{(i)}$ for~$i \in \{i_1,\dots,i_2\}$. Then, by~\eqref{eq:mul_add_exp} and the properties of the matrix exponential,
\begin{equation}
\begin{aligned}
    H_i &= \exp(A^{[k - i]}t) \otimes \exp(B^{[i]}t) \\
    &= \exp((A^{[k - i]} \oplus B^{[i]})t).
\end{aligned}
\end{equation}
Thus,
$
\frac{d}{dt}H_i |_{t=0}
 =
A^{[k-i]}  \oplus  B^{[i]},
$
and combining  this with~\eqref{eq:blk_diag_mul}
proves~\eqref{eq:blk_diag_add}.
\end{IEEEproof}

Theorem~\ref{thm:ron} allows 
to  derive a lower and an upper bound on the matrix measure of the~$k$th additive compound of a block-diagonal matrix.

\begin{Corollary}\label{cor:parallel_measure}
Let~$A\in\R^{n\times n}$ and~$B\in\R^{m\times m}$ and define $C := \diag(A,B)$. Fix~$k\in\{1,\dots,n+m\}$,
and let~$r:=\binom{n+m}{k}$. Let~$i_1:=\max\{0,k-n\} $ and~$i_2:=\min\{m,k\} $. Let $\mu_i$, $i \in \{i_1,\dots,i_2\}$, be matrix measures induced by (possibly different)~$L_{p_i}$ vector norms. Let~$|\cdot|_0$ denote some monotonic norm, define a norm~$|\cdot|$ as in~\eqref{eq:hier_norm} and let~$\mu$ denote the matrix measure induced by the norm $|\cdot|$. Let~$P\in\R^{r \times r}$ denote the permutation matrix in Thm.~\ref{thm:ron}, and define the norm $|x|_P := |P^{-1}x|$. Let~$\mu_P$ denote the matrix measure induced by $|\cdot|_P$. Then
\begin{equation}\label{eq:par_meas}
    \mu_P \left( C^{[k]} \right) = \max_{i \in \{ i_1,\dots, i_2 \}} \{ \mu_i(A^{[k-i]}) + \mu_i(B^{[i]}) \},
\end{equation}
and
\begin{equation}\label{eq:par_meas_lb}
    \min_{i \in \{ i_1,\dots, i_2 \}} \{ -\mu_i(-A^{[k-i]}) - \mu_i(-B^{[i]}) \} \le \mu_P \left( C^{[k]} \right).
\end{equation}
\end{Corollary}
\begin{IEEEproof}
By Thm.~\ref{thm:ron}
\begin{align*}
\mu_P(C^{[k]})& = \mu\left( \diag_{i\in\{i_1,\dots,i_2\}}(A^{[k-i]} \oplus B^{[i]}) \right)\\& = \max_i \mu_i(A^{[k-i]} \oplus B^{[i]}),
\end{align*}
where we used Thm.~\ref{thm:strom} and Eq.~\eqref{eq:strom_diag} for the second equality. Since $\mu_i$ are induced by $L_p$ norms, we have that (see for example Prop.~5 and Thm.~6 in~\cite{pines2021})
\begin{equation}
    \mu_i(A^{[k-i]} \oplus B^{[i]}) = \mu_i(A^{[k-i]}) + \mu_i(B^{[i]}),
\end{equation}
and this proves the equality in~\eqref{eq:par_meas}. The lower bound in~\eqref{eq:par_meas_lb} may be derived by noting that for any matrix measure and any matrix~$M$
\begin{equation}
    -\mu(-M) \le \mu(M).
\end{equation}
\end{IEEEproof}

We may now prove Thm.~\ref{thm:serial_k_contract}.
\begin{IEEEproof}[Proof of Theorem \ref{thm:serial_k_contract}]
Recall that the Jacobian~$J$ of the series connection  has the form~\eqref{eq:jaco_serial}. 
For~$\epsilon>0$,     let
$
        T(\epsilon) = \begin{bmatrix}
            I_n & 0 \\
            0 & \epsilon I_m
        \end{bmatrix}$. Then
    \[    T(\epsilon)JT^{-1}(\epsilon) = \begin{bmatrix}
            J^{11} & 0 \\
            \epsilon J^{21} & J^{22}
        \end{bmatrix}.
    \]
    Since $J^{21}(t,x)$ is uniformly bounded,
    we can make the term~$\epsilon  J ^{21}$ arbitrarily small. Let
$
        \tilde{J} := \begin{bmatrix}
            J^{11} & 0 \\
            0 & J^{22}
      \end{bmatrix} 
 $.
    By Corollary~\ref{cor:parallel_measure} and~\eqref{cond:block_contraction}, there exists a  matrix measure $\mu$ such that
    \[ \mu(\tilde{J}) \le -\min_i{\eta_i} < 0.
    \]
    Let $|\cdot|$ denote the vector norm corresponding to $\mu$,  and define a scaled vector norm by
    $ |y|_{\epsilon} := |T(\epsilon) y|$.
    Let $\mu_\epsilon$ denote the matrix measure induced by  this  scaled norm. Then for all sufficiently small $\epsilon>0$,
    \[
    \mu_\epsilon(J) \le -\min_i{\eta_i/2} < 0
    \]
    for all $x^1 \in \Omega^1, x^2 \in \Omega^2$, and $t \ge 0$, and this completes the proof.
\end{IEEEproof}

Up to this point we considered only series interconnections. 
However, for the special case of $k$-contraction w.r.t. the~$L_2$ norm it is possible to study also another form of interconnection, namely, a skew-symmetric feedback connection.

\subsection{$k$-contraction of a skew-symmetric feedback interconnection}
 Consider the system
\begin{equation}\label{eq:skew_sym_feedback}
\begin{aligned}
    \dot{x}^1 &= f^1(t, x^1, x^2), \\
    \dot{x}^2 &= f^2(t, x^1, x^2),
\end{aligned}
\end{equation}
where we assume that~$x^1$ and~$x^2$ evolve on convex sets~$\Omega^1 \subseteq \R^{n  }$ and~$\Omega^2 \subseteq \R^{m}$ respectively, and for every initial condition $a \in \Omega^1 \times \Omega^2$ a unique solution exists. Let~$J^{ij} := \frac{\partial}{\partial x^j} f^i$.
    In this section, we assume that
    there exists~$c>0$ such that
    \be\label{eq:skew_connection}
    J^{21}(t,x^1,x^2) = -c(J^{12}(t,x^1,x^2))^T
    \ee
    for all~$x^i \in \Omega^i$ and all~$t\geq 0$.
Let~$\mu_2$ denote the matrix measure induced by the~$L_2$ norm. 
\begin{Proposition}\label{thm:skew_sym_contract}
    Fix~$k \in \{1,\dots,n+m\}$, and let~$i_1 := \max\{0,k-n\}$ and~$i_2 := \min\{m,k\}$. Suppose that~\eqref{eq:skew_connection} holds. If
    \begin{equation}\label{cond:block_contraction_L2}
        \mu_2((J^{11}(t,x^1,x^2))^{[k-i]}) + \mu_2((J^{22}(t,x^1,x^2))^{[i]}) \le -\eta < 0
    \end{equation}
    for all $i \in \{i_1,\dots,i_2\}$,~$x^1 \in \Omega^1$,~$x^2 \in \Omega^2$ and~$t \ge 0$ then~\eqref{eq:skew_sym_feedback} is $k$-contracting.
\end{Proposition}

\begin{IEEEproof}
The Jacobian of~\eqref{eq:skew_sym_feedback} is
\[
J(t,x^1,x^2) = \begin{bmatrix}
        J^{11}(t,x^1,x^2) & J^{12}(t,x^1,x^2) \\
        -c(J^{12}(t,x^1,x^2))^T & J^{22}(t,x^1,x^2)
    \end{bmatrix}.
  \]
Define $T := \begin{bmatrix} \sqrt{c} I_n & 0 \\ 0 & I_m \end{bmatrix}$.
Then\[
TJT^{-1}=\begin{bmatrix}
J^{11}&\sqrt{c} J^{12} \\
-\sqrt{c}(J^{12})^T &J^{22}  
\end{bmatrix} ,
\]
and
\[
    (T J T^{-1})^{[k]} = \begin{bmatrix}
        J^{11} & 0 \\
        0 & J^{22}
    \end{bmatrix}^{[k]} + \begin{bmatrix}
        0 & \sqrt{c}J^{12} \\
        -\sqrt{c}(J^{12})^T & 0
    \end{bmatrix}^{[k]}.
\]
It follows from~\eqref{eq:matirxm_k} that~$\mu_2((TJT^{-1})^{[k]} ) =\mu_2 \left( \begin{bmatrix}
        J^{11} & 0 \\
        0 & J^{22}
    \end{bmatrix}^{[k]} \right) $,
so effectively we have   a block-diagonal Jacobian. Applying
 Corollary~\ref{cor:parallel_measure} completes the proof. 
 \end{IEEEproof}

\section{An application}\label{sec:appli}
Desoer and Haneda~\cite{desoer_1972}
have shown that 1-contracting   systems 
with an additive  input satisfy what is now known as an input-to-state stability~(ISS) property. The ISS property has  become a fundamental topic in systems and control theory~\cite{sontag_iss}. More recently,
contracting (i.e. 1-contracting) systems with inputs have   been studied in~\cite{sontag2010contractive}.

In this section,  we consider a closely related question, namely,
under what conditions a control system, with a time-varying  exponential  input, is equivalent to a time-invariant 
$k$-contracting system. This can be   studied  
using the results derived above after expressing the control system  as  the series interconnection of two sub-systems.

Consider the control  system
\begin{align}\label{eq:contsys}
\dot x&=f(x) + g(u),
\end{align}
where~$f:\R^n\to\R^n$ and $g:\R\to\R^n$ are~$C^1$.
Consider the input
\begin{align}\label{eq:input_exp}
u(t) &=\exp(\alpha t),
\end{align}
  where~$\alpha\in\R$.
  
 Write the closed-loop system
as the time-invariant $(n+1)$-dimensional system
\begin{align}
\label{eq:new_rep}
    \dot x&=f(x)+ g(y),\nonumber\\
\dot y&=\alpha y,
\end{align}
with~$y(0)=1$.
The Jacobian of~\eqref{eq:new_rep}
\[
J(x,y)= \begin{bmatrix}
        \frac{\partial}{\partial x}f(x)& \frac{\partial}{\partial y}g(y)\\
        0& \alpha
\end{bmatrix}.
\]
Applying  Theorem~\ref{thm:serial_k_contract}   yields the following result.
  \begin{Corollary}\label{core:decy_exp}
  Suppose that the trajectories of~\eqref{eq:new_rep}  evolve on a convex set~$\Omega\subset\R^{n+1}$, and that $\frac{\partial}{\partial u}g(u)$ is uniformly bounded.
Fix~$k\in\{1,\dots,n\}$. 
  If there exist~$\eta>0$ and matrix measures $\mu_i,\mu_j$   induced by (possibly different) $L_p$ norms such that
    \begin{align}
    \label{eq:two_cond} 
        \mu_i( (\frac{\partial }{\partial x}f(x) ) ^{[k]})  \le -\eta  ,\\
\mu_j( (\frac{\partial }{\partial x} f(x) )  ^{[k-1]}) + \alpha \le -\eta   ,\nonumber
    \end{align}
    for  all $x  \in \Omega $ then  the time-invariant  system~\eqref{eq:new_rep} is $k$-contracting.
  \end{Corollary}

Note that for~$k=1$, i.e. when we require
1-contractivity, this reduces to the requirement that 
the open-loop system is 1-contracting and that~$\alpha<0$, that is, an exponentially converging input. 
For~$k=2$, the first condition in~\eqref{eq:two_cond}  is that the open-loop system is $2$-contracting. If~$\alpha\geq 0 $ then the second 
condition in~\eqref{eq:two_cond}  implies that the open-loop system is also $1$-contracting. But, for~$\alpha<0$,
the second condition may hold even if the open-loop system is not $1$-contracting.

Since the system~\eqref{eq:new_rep} 
is time-invariant, specializing  Corollary~\ref{core:decy_exp}  
to the case~$k=2$ gives the following result. 
\begin{Corollary}\label{coro:exp_inp_k2}
Consider the control  system described by
\eqref{eq:contsys} and~\eqref{eq:input_exp}, and suppose that its trajectories evolve on a convex set~$\Omega \subseteq\R^n$, and that   $\frac{\partial}{\partial u}g(u)$ is uniformly bounded.
  If there exist~$\eta>0$ and matrix measures $\mu_i,\mu_j$   induced by (possibly different) $L_p$ norms such that
    \begin{align}\label{eq:k2exp} 
        \mu_i( (\frac{\partial }{\partial x}f(x) ) ^{[2]})  \le -\eta  ,\nonumber \\
\mu_j( \frac{\partial }{\partial x} f(x) )  + \alpha \le -\eta   ,
\end{align}
    for  all $x  \in \Omega $ then the system~\eqref{eq:new_rep} is~$2$-contracting,  
    and thus any bounded trajectory of the closed-loop  system converges to the set of equilibria.
\end{Corollary} 

\begin{figure}[t]
 \begin{center}
\includegraphics[width=8cm,height=7cm]{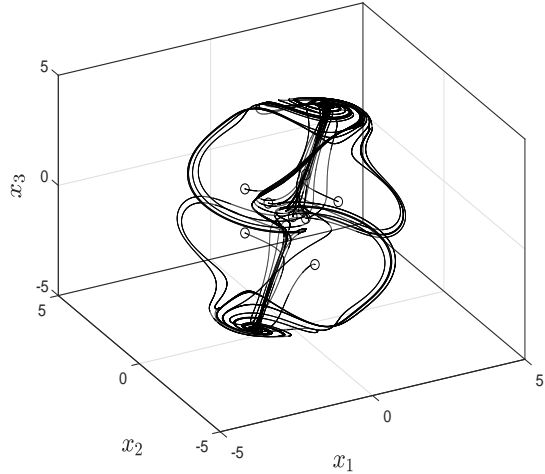}
\caption{Several trajectories of~\eqref{eq:thom} with~$d=0.193186$. The circles denote the initial conditions of the trajectories.  }\label{fig:chaos}
\end{center}
\end{figure}

\begin{figure}[t]
 \begin{center}
\includegraphics[width=8cm,height=7cm]{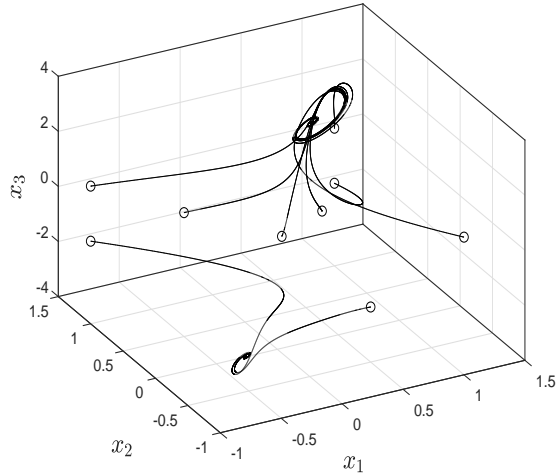}
\caption{Several trajectories of  the perturbed closed-loop      system~\eqref{eq:closedfg_withnoise}. The circles denote the initial conditions of the trajectories.  }\label{fig:chaos_closed}
\end{center}
\end{figure}

Corollary~\ref{coro:exp_inp_k2} is a kind of analogue of the result of
Desoer and Haneda~\cite{desoer_1972}
for  1-contracting   systems, as it provides conditions guaranteeing that a 2-contracting system with an exponential input retains the useful properties of time-invariant  2-contracting systems.
The next example demonstrates how these theoretical results can be used to study the robustness  of a nonlinear control system.  
\begin{Example}
Ref.~\cite{pines2021}
designed a feedback controller for 
  a popular chaotic system, introduced by Thomas~\cite{thomas99} (see also the recent review~\cite{chaos_survey}), so that the 
  closed-loop system is 
  $2$-contracting. Here we show how our results can be used to analyse the closed-loop system with an added perturbation  modeled as a decaying exponential. 
  
Thomas' cyclically symmetric attractor is given by:
\begin{align} \label{eq:thom}
\dot x_1 =&  \sin(x_2)- d x_1 , \nonumber \\
\dot x_2 =&  \sin(x_3) - d x_2, \\
\dot x_3 =& \sin(x_1) - d x_3, \nonumber
\end{align}
where $d>0$ is the
dissipation  constant. 
Note that the convex and compact set~$D: = \{x\in\R^3:  d  |x|_\infty \leq 1 \}$ is an
invariant set of the dynamics.
For~$d>1$ 
the origin is the single  equilibrium of~\eqref{eq:thom}. 
The system
undergoes a series of bifurcations as~$d$ decreases,
   and becomes chaotic at~$d  \approx 0.208186$.  
 Fig.~\ref{fig:chaos} depicts the trajectories emanating from several initial conditions for~$d=0.193186$.   It may be seen that they all converge to a  strange attractor.

 Let~$f(x)$ denote the vector field in~\eqref{eq:thom}.
 Ref.~\cite{pines2021} considered the closed-loop system
 \be\label{eq:closedfg}
 \dot x=f(x)+g(x),
 \ee
 where~$g(x)$ is the linear partial-state controller 
\[
g(x )=-\diag(c,c,0) x ,
\]
with~$c>0$.  The Jacobian  of the closed-loop system is
\[
J(x)=\begin{bmatrix}
-d-c &  \cos(x_2) &0 \\ 
0&-d-c&\cos(x_3) \\
\cos(x_1) &0&-d
\end{bmatrix},
\]
and
\[
J^{[2]}(x)=\begin{bmatrix}
-2d-2c &  \cos(x_3) &0 \\ 
0&-2d-c&\cos(x_2) \\
-\cos(x_1) &0&-2d-c
\end{bmatrix}.
\]
Thus,
\begin{align*}
\mu_1(J)&=\max \{  -d-c+|\cos(x_1)|,
-d-c+|\cos(x_2) |,\\& \qquad \;\;\;\;\;
-d+|\cos(x_3)| \} \\
&\leq -d+1,
\end{align*}
and 
\begin{align*}
\mu_1(J^{[2]})&=\max \{  -2(d+c)+|\cos(x_1)|,
-2d-c+|\cos(x_3) |,\\
&\qquad \;\;\;\;\;
-2d-c+|\cos(x_2)| \} \\
&\leq -2d-c+1.
\end{align*}
We conclude that the closed-loop system is $1$-contracting  if~$d>1$, and~$2$-contracting   if~$c>1-2d$. Note that if the closed-loop system is~$2$-contracting then in particular it admits a well-ordered behaviour, so the controller ``de-chaotifies'' the original dynamics.

Assume that the closed-loop system is perturbed by an additive exponentially decaying noise, so that the dynamics becomes
\be\label{eq:closedfg_withnoise}
 \dot x=f(x)+g(x)+b\exp(\alpha t) ,
 \ee
 with~$b\in\R^3$ and~$\alpha<0$. Since the uncontrolled system is chaotic, the perturbation may have a strong effect on the dynamical behaviour. 
 
 By Corollary~\ref{coro:exp_inp_k2},
 the perturbed system will be equivalent to a time-invariant  $2$-contracting system if
 \begin{align}\label{eq:cdcon}
 c>1-2d,\quad
 \alpha<1-d.
 \end{align}
 Fig.~\ref{fig:chaos_closed} depicts trajectories of the system~\eqref{eq:closedfg_withnoise} with~$b=(1/8)\begin{bmatrix}
 1&1&1
 \end{bmatrix}^T$,
 $d=0.193186  $, $c=1.1-2d$, and~$\alpha=-0.1$, so~\eqref{eq:cdcon} holds. It may be seen that every trajectory converges to an equilibrium point, as expected. Note that there are several equilibrium points, so the system is not~$1$-contracting  w.r.t.  any norm.
\end{Example}

\section{Discussion}
A fundamental  topic in systems theory is the analysis of  an interconnected system  based on the properties of its  sub-systems and the interconnection network.  
In this context, an
  important advantage of contracting systems is that 
various interconnections of contracting  sub-systems yield   a contracting  system.

We derived a new sufficient condition guaranteeing that the series connection of two sub-systems is~$k$-contracting. This is based on a new formula for the~$k$th compounds of a block-diagonal matrix. The latter result may find more applications. For example, recall that an LTV system  is called $k$-positive if  it maps the set of vectors with up to~$k-1$ sign variations to itself~\cite{is_my_system,Eyal_k_posi,rola_spect,alseidi2019discrete}. In particular,  $1$-positive systems are just positive systems. The conditions for~$k$-positivity 
are based on the structure of the~$k$ compounds of the system, so Thm.~\ref{thm:ron}
may perhaps be used to determine when a series connection of two sub-systems generates a $k$-positive system.

Thms.~\ref{thm:serial_k_contract} and~\ref{thm:skew_sym_contract} suggest that it might not be possible to improve the contraction order nor the contraction 
rate using a series or skew-symmetric feedback interconnection. However, these ideas may still be useful in designing a controller which retains the $k$-contraction of a system, while adding other desirable properties. For example, it might be possible to stabilize a 2-contracting system by adding a controller which guarantees that all solutions are bounded while at the same time retaining 2-contraction.

Thm.~\ref{thm:serial_k_contract} may also be used to study $k$-contraction 
in a system whose output is  fed to an integrator.
Indeed, given  $\dot{x} = f(t,x)$, let $g(x)$ denote the output of the system, and consider the augmented system
\[
\begin{aligned}
    \dot{x} &= f(t,x), \\
    \dot{y} &= g(x).
\end{aligned}
\]
 Thm.~\ref{thm:serial_k_contract} can be applied to study~$k$-contraction of  this serial interconnection.
 
As another topic for further research, recall that the theory of \emph{asymptotic autonomous systems} considers the system 
$\dot x =f(t,x) $, where 
$f( t,x) $ converges (in some technical sense) to the time-invariant vector field~$ g(x)$
as~$t\to\infty$. The question is what is the relation between the solutions~$x(t,t_0,x_0)$ and the solutions~$y(t,y_0)$ of the time-invariant system~$\dot y=g(y)$ (see e.g.~\cite{cics2006,Thieme_asym_cont_theory} and the references therein). This is   related to the question of retaining~$k$-contraction   in a system perturbed by a decaying exponential, and it might be of interest  to  further explore this connection.

\end{document}